\documentclass{article}
\usepackage{amssymb,latexsym,amscd}

\newtheorem{theorem}{Theorem}
\newtheorem{lm}{Lemma}
\newtheorem{prop}{Proposition}

\begin{document}
\hbadness=10000
\title{{\bf Deformation of integral coisotropic submanifolds in symplectic manifolds}}
\author{Wei-Dong Ruan\\
Department of Mathematics\\
University of Illinois at Chicago\\
Chicago, IL 60607\\}
\date{Revised October 2004}
\footnotetext{Partially supported by NSF Grant DMS-0104150.}
\maketitle
\begin{abstract}
In this paper we prove the unobstructedness of the deformation of integral coisotropic submanifolds in symplectic manifolds, which can be viewed as a natural generalization of Weinstein's results (\cite{W}) for Lagrangian submanifolds.\\\\
\end{abstract}

Let $(X,\omega)$ be a symplectic manifold. A submanifold $Y \subset (X,\omega)$ is coisotropic if ${\rm rank}(\omega|_Y) = 2\dim(Y) - \dim X$. Real hypersurfaces and Lagrangian submanifolds in $(X,\omega)$ are examples of coisotropic submanifolds. The closedness of $\omega$ implies that ${\rm ker}(\omega|_Y)$ defines an integrable distribution on the coisotropic submanifold $Y$. The corresponding foliation ${\cal F}$ is called null foliation. We call $Y$ an {\bf integral} coisotropic submanifold if the leaves of the null foliation ${\cal F}$ are all closed and form a fibration $\pi: Y \rightarrow S$. Denote $F_s = \pi^{-1}(s)$ for $s\in S$. Examples of integral coisotropic submanifolds are given later in this paper.\\

For a Lagrangian submanifold $L \subset (X,\omega)$, according to the work of Weinstein (\cite{W}), there exists an open neighborhood $U \subset (X,\omega)$ of $L$ and a symplectic open embedding $(U,\omega) \rightarrow (T^*L, \omega_{T^*L})$, such that $L\subset U$ is mapped to the zero section of $T^*L$, where $\omega_{T^*L}$ is the canonical symplectic form on $T^*L$. The moduli space of Lagrangian submanifolds in $(X,\omega)$ near $L$ modulo Hamiltonian deformation can be canonically identified with $H^1(L,\mathbb{R})$.\\

$(Y,\omega_Y)$ is called a pre-symplectic manifold if $\omega_Y$ is a closed 2-form on $Y$ with constant rank. Gotay's coisotropic neighborhood theorem \cite{G} implies that there exists a symplectic neighborhood $(U,\omega)$ containing $Y$ such that $\omega|_Y = \omega_Y$. Moreover, for another such symplectic neighborhood $(U',\omega')$ containing $Y$, by shrinking $U$ and $U'$ suitably, there exists a symplectomorphism $(U,\omega) \rightarrow (U',\omega')$ that fixes $Y$. A coisotropic submanifold $Y \subset (X,\omega)$ naturally gives rise to a pre-symplectic structure $(Y,\omega|_Y)$. The symplectic neighborhood of the coisotropic submanifold $Y$ is completely determined by the pre-symplectic structure $(Y,\omega|_Y)$.\\

The rather classical concept of coisotropic submanifold attracted quite some recent interests due to the work of Kapustin and Orlov \cite{KO} relating coisotropic submanifolds in Calabi-Yau manifolds to mirror symmetry. Unlike the deformation of Lagrangian submanifold modulo Hamiltonian equivalence, which is unobstructed, canonically affine and of finite dimension, the deformation of general coisotropic submanifold modulo Hamiltonian equivalence turns out to be much more non-trivial and complicated. According to the work of Oh and Park \cite{OP}, the formal deformation of general coisotropic submanifold modulo Hamiltonian equivalence can be reduced to a strongly homotopy Lie algebroid, and is usually obstructed and of infinite dimension.\\

In this work, we will discuss the deformation theory of {\bf integral} coisotropic submanifolds. We will show that the deformation of integral coisotropic submanifold modulo Hamiltonian equivalence is unobstructed, canonically affine and of finite dimension. In some sense, integral coisotropic submanifolds are more natural analogues of Lagrangian submanifolds in terms of deformation theory. We should point out that being integral coisotropic is not invariant under small deformations of coisotropic submanifolds.\\

In the appendix, we will give an alternative simple proof of the unobstructedness of the deformation of integral coisotropic submanifolds in symplectic manifolds using the smooth family version of coisotropic neighborhood theorems from \cite{G,OP}. From this point of view, the unobstructedness of the deformation of integral coisotropic submanifolds can also be viewed as a consequence of the general philosophy (discussed in \cite{OP}) that the deformation of coisotropic submanifolds is determined by the deformation of the corresponding pre-symplectic manifolds, and the unobstructedness of the deformation of symplectic manifolds.\\

The following lemma should be well known in symplectic geometry.\\
\begin{lm}
\label{be}
For an integral coisotropic submanifold $Y \subset (X,\omega)$, there exists a symplectic form $\omega_S$ on $S$ such that $\omega|_Y = \pi^*\omega_S$.
\end{lm}
{\bf Proof:} Let $W$ be a vector field on $Y$ along ${\cal F}$. Then $\imath(W)\omega|_Y =0$ and ${\cal L}_W\omega|_Y = d\imath(W)\omega|_Y + \imath(W)d\omega|_Y =0$. These give us the desired result.
\hfill\rule{2.1mm}{2.1mm}\\

Let ${\cal B}_Y$ denote the space of 1-forms $\beta$ on the integral coisotropic submanifold $Y$ such that $d\beta = \pi^* \gamma$ for certain closed 2-form $\gamma$ on $S$. Consider a smooth family $\{Y_t\}$ of integral coisotropic submanifolds in $(X,\omega)$ with $Y_0 = Y$. We choose the deformation vector fields $\{V_t\}$ for the family so that the null foliation ${\cal F}_t$ is preserved under the flow $\{\phi_t\}$ of $\{V_t\}$. (In general, such choice is possible if and only if the null foliation structure is unchanged in the family. This is obvious in our integral case, where the null foliation is a fibration.) Let $\beta_t = \imath(V_t)\omega$, and let ${\bf H}$ denote the local system on $S$ with the fibre ${\bf H}_s = H^1(F_s,\mathbb{R})$ for $s\in S$. We have\\
\begin{lm}
\label{bc}
$\beta = \beta_0$ is in ${\cal B}_Y$. Consequently, $[\beta] = \{\beta|_{F_s}\}_{s\in S}$ represents a locally constant section of ${\bf H}$ over $S$. Conversely, any element of $H^0(S,{\bf H})$ can be represented in such way by some $\beta \in {\cal B}_Y$ so that $[\beta] \rightarrow [\gamma]$ defines the natural map $H^0(S,{\bf H}) \rightarrow H^2(S,\mathbb{R})$. $[\beta] = 0$ in $H^0(S,{\bf H})$ if and only if $\beta$ determines an infinitesimal Hamiltonian deformation.
\end{lm}
{\bf Proof:} Let $W$ be a vector field on $Y$ along ${\cal F}$. Since the null foliation is preserved by the flow, we have $\imath(W)\phi_t^*\omega|_Y=0$, $\imath(W)d\beta = \imath(W){\cal L}_{V}\omega|_Y=0$. Consequently, $d\beta$ descends to a closed 2-form $\gamma$ on $S$.\\

Locally on $S$, we may write $\gamma = d\tilde{\beta}$. Then $\beta - \pi^*\tilde{\beta}$ is closed and $\beta|_{F_s} = (\beta - \pi^*\tilde{\beta})|_{F_s}$ represents a locally constant class in $H^1(F_s,\mathbb{R})$ for $s\in S$.\\

For the converse, pick a covering $\{U_\alpha\}$ for $S$ by open sets whose intersections are contractible, with a compatible partition of unity $\{\rho_\alpha\}$. For $b \in H^0(S,{\bf H})$, it is easy to find closed 1-form $\beta'_\alpha$ on $\pi^{-1}(U_\alpha)$ such that $\beta'_\alpha|_{F_s}$ represents $b_s \in H^1(F_s,\mathbb{R})$ for $s\in U_\alpha$. $\beta' = \sum_\alpha \beta'_\alpha\pi^*\rho_\alpha$ defines a global 1-form on $Y$ such that $\beta'|_{F_s}$ represents $b_s \in H^1(F_s,\mathbb{R})$ for $s\in S$. Consequently, there exists smooth function $f_\alpha$ on $\pi^{-1}(U_\alpha)$ such that $(\beta' - \beta'_\alpha)|_{F_s} = df_\alpha|_{F_s}$ for $s\in U_\alpha$. Let $\beta_\alpha = \beta'_\alpha + df_\alpha$, then $\beta_\alpha$ is a closed 1-form on $\pi^{-1}(U_\alpha)$ such that $\beta_\alpha|_{F_s} = \beta'|_{F_s}$ represents $b_s \in H^1(F_s,\mathbb{R})$ for $s\in U_\alpha$. Hence, $\beta_{\alpha\alpha'} = \beta_\alpha - \beta_{\alpha'}$ defined over $\pi^{-1}(U_\alpha \cap U_{\alpha'})$ is closed and vanishes when restricted to $F_s$ for $s\in U_\alpha \cap U_{\alpha'}$. Consequently, there exists a closed 1-form $\hat{\beta}_{\alpha\alpha'}$ on $U_\alpha \cap U_{\alpha'}$ such that $\beta_{\alpha\alpha'} = \pi^*\hat{\beta}_{\alpha\alpha'}$. Since the sheaf of 1-forms on $S$ is soft, there exists 1-form $\hat{\beta}_{\alpha}$ on each $U_\alpha$ such that $\hat{\beta}_{\alpha\alpha'} = \hat{\beta}_\alpha - \hat{\beta}_{\alpha'}$. By our construction, it is easy to see that $\beta_\alpha - \pi^* \hat{\beta}_\alpha$ on $\pi^{-1}(U_\alpha)$ (resp. $-d\hat{\beta}_\alpha$ on $U_\alpha$) patch up to a global 1-form $\beta$ on $Y$ (resp. global closed 2-form $\gamma$ on $S$) so that $d\beta = \pi^*\gamma$.\\

$[\beta] = 0$ in $H^0(S,{\bf H})$ if and only if $\beta|_{F_s}$ is exact for all $s\in S$. Hence, there exists a smooth function $f$ on $Y$ such that $\beta|_{F_s} = df|_{F_s}$ for all $s\in S$. Let $W$ be a vector field on $Y$ along ${\cal F}$. Then $\imath(W) (\beta - df) =0$ and ${\cal L}_W (\beta - df) = d\imath(W)(\beta - df) + \imath(W)d(\beta - df) =\imath(W)\pi^*\gamma =0$. Consequently, there exists a 1-form $\breve{\beta}$ on $S$ such that $\beta - df = \pi^* \breve{\beta}$, $\beta = df + \pi^* \breve{\beta}$. Since the vector field determined by $\pi^* \breve{\beta}$ is tangent to $Y$, $\beta$ determines an infinitesimal Hamiltonian deformation.
\hfill\rule{2.1mm}{2.1mm}\\

Lemma \ref{bc} implies that the infinitesimal deformations of an integral coisotropic submanifold $Y$ modulo infinitesimal Hamiltonian equivalence lie in $H^0(S,{\bf H})$. In the following proposition, we will verify that each $[\beta] \in H^0(S,{\bf H})$ can indeed be realized as the infinitesimal deformation of a family $\{Y_t\}$ of integral coisotropic submanifolds with $Y_0=Y$. Namely, the deformation of integral coisotropic submanifolds is unobstructed.\\
\begin{prop}
\label{bd}
For any $[\beta]\in H^0(S,{\bf H})$, there exists a family $\{Y_t\}$ of integral coisotropic submanifolds for small $t$ with $Y_0=Y$, whose infinitesimal deformation at all $t$ is represented by $[\beta]$.
\end{prop}
{\bf Proof:} For any $U_\alpha$ in the open covering $\{U_\alpha\}$ of $S$, assume the deformation family $\{Y_t\}$ of integral coisotropic submanifolds with $Y_0 = \pi^{-1}(U) \subset Y$ is constructed for some open subset $U\subset U_\alpha$ such that the infinitesimal deformation at any $t$ is represented by $[\beta]|_{U} \in H^0(U,{\bf H})$. In the following, we will extend the deformation family to $\{Y_t\}$ with $Y_0 = \pi^{-1}(U_\alpha)$.\\

According to lemma \ref{bc}, we may assume $d\beta = \pi^*\gamma$ for some closed 2-form $\gamma$ on $S$. Since $U_\alpha$ is contractible, there exists 1-form $\beta'_\alpha$ on $U_\alpha$ so that $\gamma |_{U_\alpha} = d\beta'_\alpha$. Then the closed 1-form $\beta_\alpha = \beta|_{\pi^{-1} (U_\alpha)} - \pi^* \beta'_\alpha$ on $\pi^{-1}(U_\alpha)$ satisfies $\beta_\alpha|_{F_s} = \beta|_{F_s}$ for $s\in U_\alpha$.\\

Choose the deformation vector fields $\{V_t\}$ for the family $\{Y_t\}$ so that the null foliation ${\cal F}_t$ is preserved under the flow $\{\phi_t\}$ of $\{V_t\}$. The flow gives the smooth family of identifications $\phi_t: Y_0 = \pi^{-1}(U) \rightarrow Y_t$ for $t$ small that also identify the fibres of null foliation. Let $\beta_t = \phi_t^* \imath (V_t) \omega$. By our assumption, $[\beta_t] = [\beta]|_{U} = [\beta_\alpha]|_{U} \in H^0(U,{\bf H})$. It is straightforward to see that this condition implies that there exists a smooth family of functions $\{f_t\}$ on $Y_t \cong \pi^{-1}(U)$ and 1-forms $\beta'_t$ on $U$ such that $\beta_t = \beta_\alpha|_{\pi^{-1}(U)} + df_t + \pi^* \beta'_t$ and $f_0=0$. Replacing $\beta_t$ by $\beta_\alpha|_{\pi^{-1}(U)} + df_t$, we can make $\beta_t$ closed and $\beta_0 = \beta_\alpha|_{\pi^{-1}(U)}$. $\{V_t\}$ still generates the same $\{Y_t\}$ because the vector field corresponding to $\pi^* \beta'_t$ is tangent to $Y_t$.\\

Let $\hat{U}_\alpha$ be a small neighborhood of $\pi^{-1}(U_\alpha)$ in $X$. Since $\{Y_t\}$ is very close to $Y_0 = \pi^{-1}(U)$ for $t$ small, we may construct projections $\hat{\pi}_t: \hat{U}_\alpha \rightarrow \pi^{-1}(U_\alpha)$ so that $\hat{\pi}_t|_{Y_t} = \phi_t^{-1}: Y_t \rightarrow Y_0 =\pi^{-1}(U)$ and $\hat{\pi}_0|_{\pi^{-1}(U_\alpha)} = {\rm id}_{\pi^{-1}(U_\alpha)}$. (Some of the open sets may need to be shrunk slightly, which does not affect our proof.)\\

Extending $\{f_t\}$ from $\pi^{-1}(U)$ to a smooth family of functions on $\pi^{-1}(U_\alpha)$ still satisfies $f_0=0$. Then $\hat{\beta}_t = \hat{\pi}_t^* (\beta_\alpha + df_t)$ is a smooth family of closed 1-forms on $\hat{U}_\alpha$ satisfying $\hat{\beta}_t|_{Y_t} = \imath (V_t) \omega$ and $\hat{\beta}_0|_{\pi^{-1}(U_\alpha)} = \beta_\alpha$ for small $t$. Then $\hat{V}_t$ satisfying $\imath(\hat{V}_t)\omega = \hat{\beta}_t$ gives rise to a flow of symplectomorphisms on $\hat{U}_\alpha$. This flow generates the deformation family $\{Y_t\}$ with $Y_0 = \pi^{-1}(U_\alpha)$ that extends the previous deformation family $\{Y_t\}$ with $Y_0 = \pi^{-1}(U)$. From the expression of $\hat{\beta}_t$, the infinitesimal deformation at any $t$ is clearly represented by $[\beta]|_{U_\alpha} \in H^0(U_\alpha,{\bf H})$.\\

With the above argument, we may start the induction with $U=\emptyset$. After the deformation family $\{Y_t\}$ is constructed over $\pi^{-1}(\tilde{U}_A)$, where $\tilde{U}_A = \displaystyle\bigcup_{\alpha'\in A}U_{\alpha'}$, take $U = \tilde{U}_A \cap U_\alpha$ for $\alpha \not\in A$. By induction, we can construct the family $\{Y_t\}$ globally.
\hfill\rule{2.1mm}{2.1mm}\\

Let ${\cal M}_Y$ denote the local moduli space of integral coisotropic submanifolds (that are small deformations of $Y$) modulo Hamiltonian equivalence. To prove the smoothness of ${\cal M}_Y$ near $Y$, we need the following Proposition \ref{bf}.\\

Consider two 1-parameter families $\{Y_{r,t}\}$ of integral coisotropic submanifolds with $Y_{r,0} = Y$ for $r=0,1$. Assume the deformation vector fields $\{V_{r,t}\}$ determine the flows $\phi_{r,t}: Y \rightarrow Y_{r,t}$ that preserve the null foliations. Then the corresponding 1-forms $\beta_{r,t} = \phi_{r,t}^* \imath(V_{r,t}) \omega$ on $Y$ is in ${\cal B}_Y$. Assume that there is a smooth family of 1-forms $\{\beta_{r,t}\} \subset {\cal B}_Y$ for $r\in [0,1]$ and $t\in \mathbb{R}$ small that extends $\{\beta_{r,t}\}$ for $r=0,1$. We have\\
\begin{prop}
\label{bf}
There exists a smooth family of integral coisotropic submanifolds $\{Y_{r,t}\}$ for $r\in [0,1]$ and $t\in \mathbb{R}$ small satisfying $Y_{r,0} = Y$ that extend $\{Y_{r,t}\}$ for $r=0,1$, such that the infinitesimal deformation along $t$ direction at $Y_{r,t}$ is represented by $[\beta_{r,t}] \in H^0(S,{\bf H})$.
\end{prop}
{\bf Proof:} The proof of this proposition uses essentially the same idea as that of Proposition \ref{bd}. For any $U_\alpha$ in the open covering $\{U_\alpha\}$ of $S$, assume the deformation family $\{Y_{r,t}\}$ of integral coisotropic submanifolds with $Y_{r,0} = \pi^{-1}(U) \subset Y$ (resp. $Y_{r,0} = Y$) for $r \in (0,1)$ (resp. $r= 0,1$) extending $\{Y_{r,t}\}$ for $r=0,1$ is constructed for some open subset $U\subset U_\alpha$ such that the infinitesimal deformation at $Y_{r,t}$ for $r \in (0,1)$ is represented by $[\beta_{r,t}]|_{U} \in H^0(U,{\bf H})$. More precisely, there exist deformation vector fields $\{V_{U,r,t}\}$ and the flow $\{\phi_{U,r,t}\}$ with $\beta_{U,r,t} = \phi_{U,r,t}^* \imath (V_{U,r,t}) \omega \in {\cal B}_{\pi^{-1}(U)}$ such that $\beta_{U,r,t} = \beta_{r,t}|_{\pi^{-1} (U)}$ for $r=0,1$ and $[\beta_{U,r,t}] = [\beta_{r,t}]|_{U} \in H^0(U,{\bf H})$ for $r \in (0,1)$. It is straightforward to see that these conditions imply that there exists a smooth family of functions $\{f_{U,r,t}\}$ on $Y_t \cong \pi^{-1}(U)$ and 1-forms $\beta'_{U,r,t}$ on $U$ such that $\beta_{U,r,t} = \beta_{r,t}|_{\pi^{-1}(U)} + df_{U,r,t} + \pi^* \beta'_{U,r,t}$ and $f_{U,r,0}=f_{U,0,t}=f_{U,1,t} = 0$. In the following, we will extend the deformation family to $\{Y_{r,t}\}$ with $Y_{r,0} = \pi^{-1}(U_\alpha)$ for $r \in (0,1)$.\\

$\beta_{r,t} \in {\cal B}_Y$ implies that $d\beta_{r,t} = \pi^*\gamma_{r,t}$ for some closed 2-form $\gamma_{r,t}$ on $S$. Since $U_\alpha$ is contractible, there exists 1-form $\beta'_{\alpha,r,t}$ on $U_\alpha$ so that $\gamma_{r,t} |_{U_\alpha} = d\beta'_{\alpha,r,t}$. Then the closed 1-form $\beta_{\alpha,r,t} = \beta_{r,t}|_{\pi^{-1} (U_\alpha)} - \pi^* \beta'_{\alpha,r,t}$ on $\pi^{-1}(U_\alpha)$ satisfies $\beta_{\alpha,r,t}|_{F_s} = \beta_{r,t}|_{F_s}$ for $s\in U_\alpha$. Replacing $\beta_{U,r,t}$ by $\beta_{\alpha,r,t}|_{\pi^{-1}(U)} + df_{U,r,t} = \beta_{U,r,t} - \pi^* (\beta'_{\alpha,r,t}|_U + \beta'_{U,r,t})$, we can make $\beta_{U,r,t}$ closed and $\beta_{U,r,0} = \beta_{\alpha,r,0}|_{\pi^{-1}(U)}$. $\{V_{U,r,t}\}$ still generates the same $\{Y_{r,t}\}$ because the vector field corresponding to $\pi^* (\beta'_{\alpha,r,t}|_U + \beta'_{U,r,t})$ is tangent to $Y_{r,t}$.\\

Let $\hat{U}_\alpha$ be a small neighborhood of $\pi^{-1}(U_\alpha)$ in $X$. Since $\{Y_{r,t}\}$ is very close to $Y_{r,0} = \pi^{-1}(U)$ for $t$ small, we may construct projections $\hat{\pi}_{r,t}: \hat{U}_\alpha \rightarrow \pi^{-1}(U_\alpha)$ so that $\hat{\pi}_{r,t}|_{Y_{r,t}} = \phi_{U,r,t}^{-1}: Y_{r,t} \rightarrow Y_{r,0} =\pi^{-1}(U)$ (resp. $\hat{\pi}_{r,t}|_{Y_{r,t} \cap \pi^{-1}(U_\alpha)} = \phi_{r,t}^{-1}|_{Y_{r,t} \cap \pi^{-1}(U_\alpha)}: Y_{r,t} \cap \pi^{-1}(U_\alpha) \rightarrow Y_{r,0} \cap \pi^{-1}(U_\alpha) = \pi^{-1}(U_\alpha)$) for $r\in (0,1)$ (resp. $r=0,1$) and $\hat{\pi}_{r,0}|_{\pi^{-1}(U_\alpha)} = {\rm id}_{\pi^{-1}(U_\alpha)}$ for $r\in [0,1]$. (Some of the open sets may need to be shrunk slightly, which does not affect our proof.)\\

Extending $\{f_{U,r,t}\}$ from $\pi^{-1}(U)$ to a smooth family of functions $\{f_{\alpha,r,t}\}$ on $\pi^{-1}(U_\alpha)$ still satisfy $f_{\alpha,r,0}=f_{\alpha,0,t}=f_{\alpha,1,t} = 0$. Then $\hat{\beta}_{\alpha,r,t} = \hat{\pi}_{r,t}^* (\beta_{\alpha,r,t} + df_{\alpha,r,t})$ is a smooth family of closed 1-form on $\hat{U}_\alpha$ satisfying $\hat{\beta}_{\alpha,r,0}|_{\pi^{-1}(U_\alpha)} = \beta_{\alpha,r,0}$ and $\phi_{U,r,t}^* (\hat{\beta}_{\alpha,r,t}|_{Y_{r,t}}) = \beta_{U,r,t}$ (resp. $\phi_{r,t}^* (\hat{\beta}_{\alpha,r,t}|_{Y_{r,t} \cap \pi^{-1}(U_\alpha)}) = \beta_{r,t}|_{\pi^{-1}(U_\alpha)}$) for $r\in (0,1)$ (resp. $r=0,1$) and small $t$. Then $\hat{V}_{\alpha,r,t}$ satisfying $\imath(\hat{V}_{\alpha,r,t})\omega = \hat{\beta}_{\alpha,r,t}$ give rise to flows of symplectomorphisms on $\hat{U}_\alpha$. These flows generate the deformation family $\{Y_{r,t}\}$ with $Y_{r,0} = \pi^{-1}(U_\alpha)$ that extends the previous deformation family $\{Y_{r,t}\}$ with $Y_{r,0} = \pi^{-1}(U)$. From the expression of $\hat{\beta}_{\alpha,r,t}$, clearly, the infinitesimal deformation along $t$ direction at $Y_{r,t}$ is represented by $[\beta_{r,t}]|_{U_\alpha} \in H^0(U_\alpha,{\bf H})$.\\

With the above argument, we may start the induction with $U=\emptyset$. After the deformation family $\{Y_{r,t}\}$ is constructed over $\pi^{-1}(\tilde{U}_A)$, where $\tilde{U}_A = \displaystyle\bigcup_{\alpha'\in A}U_{\alpha'}$, take $U = \tilde{U}_A \cap U_\alpha$ for $\alpha \not\in A$. By induction, we can construct the family $\{Y_{r,t}\}$ globally.
\hfill\rule{2.1mm}{2.1mm}\\

A submanifold $F \subset (X,\omega)$ is isotropic if $\omega|_F = 0$. Let ${\cal M}_F$ denote the local moduli space of isotropic submanifolds near $F$ modulo Hamiltonian equivalence. The following proposition should be well known.\\
\begin{prop}
\label{ba}
For an isotropic submanifold $F\subset (X,\omega)$, there exists a canonical local open embedding $i_F: {\cal M}_F \rightarrow H^1(F,\mathbb{R})$.
\end{prop}
{\bf Proof:} Let $\{F_t\}$ be a deformation family of isotropy submanifolds with $F_0=F$ and $\{V_t\}$ be a smooth family of vector fields that generate the family, where $V_t$ is only defined on $F_t$. Then $\omega|_{F_t}=0$, $d\imath(V_t)\omega|_{F_t} = {\cal L}_{V_t}\omega|_{F_t}=0$.\\

Let $U$ be a small tubular neighborhood of $F$. Then it is straightforward to see that for $t$ small, each $V_t$ can be extended to a vector field $V'_t$ on $U$ so that $\imath(V'_t)\omega$ is closed and $\{V'_t\}$ forms a smooth family. Consequently, the family of isotropic submanifolds $\{F_t\}$ is generated by the symplectic flow on $U$ generated by $\{V'_t\}$. Since $H^1(U,\mathbb{R}) \cong H^1(F,\mathbb{R})$, a closed 1-form on $U$ that restrict to zero on $F$ is an exact 1-form. Hence a symplectic flow on $U$ fixing $F$ will be a Hamiltonian flow. Consequently, ${\cal M}_F$ can be identified with the local moduli space of symplectic automorphisms of $U$ near the identity map modulo Hamiltonian automorphisms of $U$, which can be identified with the local moduli space of Lagrangian submanifolds of $(U,\omega|_U)\times (U,-\omega|_U)$ near diagonal modulo Hamiltonian equivalence, which can be canonically identified with $H^1(U,\mathbb{R}) \cong H^1(F,\mathbb{R})$ locally according to \cite{W}.
\hfill\rule{2.1mm}{2.1mm}\\
\begin{theorem}
\label{bb}
There exists a canonical local open embedding $i_Y: {\cal M}_Y \rightarrow H^0(S,{\bf H})$, under which the smooth orbit of symplectomorphisms can be locally identified with the kernel of $H^0(S,{\bf H}) \rightarrow H^2(S,\mathbb{R})$.
\end{theorem}
{\bf Proof:} Notice that any two isotropic fibres in an integral coisotropic submanifold $Y$ are Hamiltonian equivalent in $X$. Take $F$ to be an isotropic fibre of $Y$. There is the natural map $\mu: {\cal M}_Y \rightarrow {\cal M}_F$. For any $[Y'] \in {\cal M}_Y$, $\mu([Y']) = [F']$, where $F'$ is an isotropic fibre in $Y'$. Proposition \ref{bd} implies that locally $i_F \circ \mu$ maps ${\cal M}_Y$ surjectively onto $H^0(S,{\bf H}) \hookrightarrow H^2(S,\mathbb{R})$. To prove the first part of the theorem, we only need to show that $i_F \circ \mu$ is injective. More precisely, for two 1-parameter families $\{Y_{r,t}\}_{r=0,1}$ of integral coisotropic submanifolds with $Y_{0,0} = Y_{1,0} = Y$, if $i_F \circ \mu(Y_{0,t_0}) = i_F \circ \mu(Y_{1,t_0})$, then $Y_{0,t_0}$ is Hamiltonian equivalent to $Y_{1,t_0}$. Assume the deformation 1-forms for the 2 families are $\{\beta_{r,t}\}_{r=0,1}$.\\

Consider the natural projection $\pi_Y: {\cal B}_Y \rightarrow H^0(S,{\bf H})$. $c_r(t) = i_F \circ \mu(Y_{r,t})$ for $r=0,1$ form 2 paths in $H^0(S,{\bf H})$ with $\frac{d c_r(t)}{dt} = [\beta_{r,t}] \in H^0(S,{\bf H})$. For $r=0,1$, $\tilde{c}_r(t)$ satisfying $\tilde{c}_r(0) = 0$, $\frac{d \tilde{c}_r(t)}{dt} = \beta_{r,t}$ is a path in ${\cal B}_Y$ that lifts $c_r(t)$, namely $\pi_Y(\tilde{c}_r(t)) = c_r(t)$. It is straightforward to construct a smooth family of paths $\{c_r(t)\}_{r\in [0,1]}$ in $H^0(S,{\bf H})$ that extend $\{c_r(t)\}_{r=0,1}$ such that $c_r(0) =0$ and $c_r(t_0) = c_0(t_0)$ for $r\in [0,1]$. Since $\pi_Y$ is linear map, one can lift $\{c_r(t)\}_{r\in [0,1]}$ to a smooth family of paths $\{\tilde{c}_r(t)\}_{r\in [0,1]}$ in ${\cal B}_Y$ extending $\{\tilde{c}_r(t)\}_{r=0,1}$ such that $\tilde{c}_r(0) =0$ and $\pi_Y(\tilde{c}_r(t)) = c_r(t)$ for $r\in [0,1]$. Let $\beta_{r,t} = \frac{d \tilde{c}_r(t)}{dt}$, then the family of 1-forms $\{\beta_{r,t}\} \subset {\cal B}_Y$ for $r\in [0,1]$ and $t\in \mathbb{R}$ small extends $\{\beta_{r,t}\}$ for $r=0,1$.\\

According to proposition \ref{bf}, there exists a smooth family of integral coisotropic submanifolds $\{Y_{r,t}\}$ for $r\in [0,1]$ and $t\in \mathbb{R}$ small satisfying $Y_{r,0} = Y$ that extend $\{Y_{r,t}\}$ for $r=0,1$, such that the infinitesimal deformation along $t$ direction at $Y_{r,t}$ is represented by $[\beta_{r,t}] \in H^0(S,{\bf H})$. Consequently, $i_F \circ \mu(Y_{r,t}) = c_r(t)$ for all $r$ and $t$. In particular, $i_F \circ \mu(Y_{r,t_0}) = c_0(t_0)$ for $r\in [0,1]$, the infinitesimal deformation along $r$ direction represents $0 \in H^0(S,{\bf H})$ and is Hamiltonian by lemma \ref{bc}. Therefore $Y_{0,t_0}$ is Hamiltonian equivalent to $Y_{1,t_0}$. This proves the first part of the theorem. The second part of the theorem is quite obvious in light of lemma \ref{bc}.
\hfill\rule{2.1mm}{2.1mm}\\

{\bf Example:} One of the simplest non-trivial integral coisotropic submanifold is the unit sphere $S^3$ in $\mathbb{R}^4$ with the standard symplectic form. The Hopf fibration gives the null foliation, with the great circles being the isotropic fibres. The symplectic quotient $S \cong S^2 = S^3/S^1$. ${\bf H}$ is trivial over $S^2$ with fibres identified with $H^1(S^1, \mathbb{R}) \cong \mathbb{R}$. $H^0(S, {\bf H}) \cong \mathbb{R}$, $H^0(S, {\bf H}) \rightarrow H^2(S, \mathbb{R})$ is an isomorphism. The moduli space ${\cal M}_{S^3} \cong \mathbb{R}_{>0}$ parameterizes round spheres of radius $r$ for $r >0$. The map ${\cal M}_{S^3} \rightarrow H^0(S, {\bf H})$ that maps the unit sphere to 0 is $r \rightarrow \frac{1}{2}(r^2 -1)$. More generally, one can consider a symplectic manifold with a symplectic real torus action. The fibres of the moment map are natural examples of integral coisotropic submanifolds.\\\\

{\large \bf Appendix: Alternative proof of proposition \ref{bd}}\\

We will start with two propositions concerning the deformation of general coisotropic submanifolds (resp. pre-symplectic manifolds) that are not necessarily integral.\\
\begin{prop}
\label{ca}
Suppose that there exists a smooth (or $C^1$) family of coisotropic
submanifolds $\{Y_t\}$ in $(X,\omega)$. $\{Y_t\}$ is generated by local symplectomorphisms if and only if $(Y_t,\omega|_{Y_t})$ are equivalent as presymplectic manifolds.
\end{prop}
{\bf Proof:} Assume that $(Y_t,\omega|_{Y_t})$ are equivalent as pre-symplectic manifolds. Then there exists a family of diffeomorphisms $\phi_t: Y_0 \rightarrow Y_t$ such that $\phi_t^*\omega|_{Y_t}= \omega|_{Y_0}$. Let $\{V_t\}$ be a smooth family of vector fields that generate the flow $\phi_t$, where $V_t$ is only defined on $Y_t$. Then $d\imath(V_t)\omega|_{Y_t} = {\cal L}_{V_t}\omega|_{Y_t}=0$.\\

Let $U$ be a small tubular neighborhood of $Y_0$. Then it is straightforward to see that for $t$ small, each $V_t$ can be extended to a vector field $V'_t$ on $U$ so that $\imath(V'_t)\omega$ is closed and $\{V'_t\}$ forms a smooth family. Consequently, the family of coisotropic submanifolds $\{Y_t\}$ is generated by the local symplectic flow on $U$ generated by $\{V'_t\}$. The other direction of the proposition is obvious.
\hfill\rule{2.1mm}{2.1mm}\\
\begin{prop}
\label{cb}
Let $Y\subset (X,\omega)$ be a coisotropic submanifold, and $\{(Y_t, \omega_t)\}$ (for $t$ in a small neighborhood of $0\in \mathbb{R}^l$) be a family of pre-symplectic manifolds such that $Y_0 = Y$, $\omega_0 = \omega|_Y$ and $[\omega_t] \in H^2(Y_t) \cong H^2(Y)$ is constant class. Then there exists a family of coisotropic embeddings $Y_t\rightarrow (X,\omega)$ for small $t$, such that $\omega_t = \omega|_{Y_t}$.
\end{prop}
{\bf Proof:} Let $U$ be a small tubular neighborhood of $Y$ in $X$. According to \cite{G,OP}, there exist a smooth family of coisotropic embeddings $\phi_t: Y_t \rightarrow (U,\tilde{\omega}_t)$ such that $\tilde{\omega}_t|_{Y_t} = \omega_t$ and $\tilde{\omega}_0 = \omega|_U$. By assumption, $[\tilde{\omega}_t] \in H^2(U) \cong H^2(Y)$ represent constant class. There exists a smooth family of symplectic morphisms $\psi_t: (U,\tilde{\omega}_t) \rightarrow (X,\omega)$ with $\psi_0 = {\rm id}_U$. $\psi_t\circ\phi_t$ will give us the desired family of coisotropic embeddings for small $t$.
\hfill\rule{2.1mm}{2.1mm}\\

{\bf Alternative proof of Proposition \ref{bd}:} Let $\rho_1: H^0(S,{\bf H}) \rightarrow H^2(S,\mathbb{R})$ be the natural map. The integral deformation of integral pre-symplectic manifold $(Y,\omega|_Y)$ is equivalent to the deformation of symplectic manifold $(S,\omega_S)$, which is well known to be unobstructed and locally parameterized by $H^2(S,\mathbb{R})$. Hence there exists a family of pre-symplectic manifolds $\{(Y_t, \omega_t)\}$ such that $Y_0 = Y$, $\omega_0 = \omega|_Y$ and the infinitesimal deformation is represented by $\rho_1([\beta])$. Since the image of $\rho_1$ is in the kernel of $\rho_2: H^2(S,\mathbb{R})\rightarrow H^2(Y,\mathbb{R})$, proposition \ref{cb} implies that there exists a family of coisotropic embeddings $Y_t\rightarrow (X,\omega)$ for small $t$, such that $\omega_t = \omega|_{Y_t}$. Let the infinitesimal deformation of the integral coisotropic sybmanifolds be represented by $[\beta'_t]$ at $Y_t$. Then $\rho_1([\beta]) = \rho_1([\beta'_t])$ for all $t$. By lemma \ref{bc}, we may assume $\beta - \beta'_t$ to be closed for all $t$. Consequently, by adjusting the family of integral coisotropic submanifolds $\{Y_t\}$ with local symplectic morphisms defined on a neighborhood of $Y_t$, one can ensure that the infinitesimal deformation at $Y_t$ be represented by $[\beta]$ for all $t$.
\hfill\rule{2.1mm}{2.1mm}\\

{\bf Remark:} This simple proof is based on suggestions from Prof. Yong-Geun Oh. The major difference of this proof from the original proof is that the original proof does not use the coisotropic neighborhood theorems from \cite{G,OP}, while this proof need the smooth family version of coisotropic neighborhood theorems from \cite{G,OP}. Of course, it is also interesting to see if Proposition \ref{bf} can be proved by propositions \ref{ca} and \ref{cb}. Although Proposition \ref{cb} was not explicitly formulated and proved in \cite{OP}, the essential idea behind Proposition \ref{cb} in some other form is already mentioned in section 8 of \cite{OP}, which modulo some standard arguments in symplectic geometry should imply Proposition \ref{cb}.\\

{\bf Acknowledgement:} This work was inspired by a talk given by Prof. Yong-Geun Oh on \cite{OP}. I would like to thank Prof. Oh for helpful discussions. I would also like to thank the referees for valuable suggestions.\\

\ifx\undefined\bysame
\newcommand{\bysame}{\leavevmode\hbox to3em{\hrulefill}\,}
\fi

\end{document}